\RequirePackage{ifpdf}
\ifpdf 
\documentclass[pdftex]{sigma}
\else
\documentclass{sigma}
\fi

\newcommand\R{{\mathbb R}}
\newcommand\RP{\R{\rm P}}

\newcommand\Sp{{\rm Sp}}
\newcommand\PSL{{\rm PSL}}
\newcommand\SL{{\rm SL}}
\newcommand\gl{{\rm gl}}

\newcommand\nor[1]{|\!|#1|\!|^2}

\newcommand\inv{k}

\newcommand\g{{\mathfrak g}}
\newcommand\n{{\mathfrak n}}
\newcommand\h{{\mathfrak h}}
\newcommand\sla{{\mathfrak sl}}
\newcommand\gla{{\mathfrak gl}}
\newcommand\oa{{\mathfrak o}}
\newcommand\spa{{\mathfrak sp}}

\newcommand\F{{\mathcal F}}
\newcommand\M{{\mathcal M}}
\newcommand\N{{\mathcal N}}
\newcommand\s{{\mathcal S}}
\newcommand\K{{\mathcal K}}
\newcommand\Pm{{\mathcal P}}
\newcommand\calD{{\mathcal D}}

\newcommand\Lo{{\mathcal L}}
\newcommand\D{{\mathcal D}}
\newcommand\HH{{\mathcal H}}

\newcommand\Z{{\mathbb Z}}

\newcommand\semip{\ltimes}

\newcommand\B{{\mathcal B}}

\newcommand\var[1]{\frac{\delta #1}{\delta L}(L)}
\newcommand\hk{\frac{\delta h}{\delta \kappa}}
\newcommand\hta{\frac{\delta h}{\delta \tau}}
\newcommand\fk{\frac{\delta f}{\delta \kappa}}
\newcommand\ft{\frac{\delta f}{\delta \tau}}

\newcommand\rr{{\bf r}}

\begin{document}
\allowdisplaybreaks

\renewcommand{\PaperNumber}{034}

\FirstPageHeading

\renewcommand{\thefootnote}{$\star$}

\ShortArticleName{Geometric Realizations of Bi-Hamiltonian Integrable Systems}

\ArticleName{Geometric Realizations of Bi-Hamiltonian\\ Completely Integrable Systems\footnote{This paper is a contribution to the Proceedings
of the Seventh International Conference ``Symmetry in Nonlinear
Mathematical Physics'' (June 24--30, 2007, Kyiv, Ukraine). The
full collection is available at
\href{http://www.emis.de/journals/SIGMA/symmetry2007.html}{http://www.emis.de/journals/SIGMA/symmetry2007.html}}}

\Author{Gloria MAR\'I~BEFFA}

\AuthorNameForHeading{G. Mar\'\i~Bef\/fa}

\Address{Department of Mathematics, University of Wisconsin, Madison, WI 53705, USA}
\Email{\href{mailto:maribeff@math.wisc.edu}{maribeff@math.wisc.edu}}
\URLaddress{\url{http://www.math.wisc.edu/~maribeff/}}

\ArticleDates{Received November 14, 2007, in f\/inal form March
13, 2008; Published online March 27, 2008}

\Abstract{In this paper we present an overview of the connection between completely in\-tegrable systems
and the background geometry of the f\/low. This relation is better seen when using a group-based concept of moving frame introduced by Fels and Olver in [{\it Acta Appl. Math.} {\bf 51} (1998), 161--213;
{\bf 55} (1999), 127--208]. The paper discusses the close connection
between dif\/ferent types of geometries and the type of equations they realize. In particular, we describe the direct relation between symmetric spaces and equations of
KdV-type, and the possible geometric origins of this connection.}


\Classification{37K25; 53A55}

\section{Introduction}

\begin{example}One of the simplest examples of a geometric realization of a completely integrable system is that of the nonlinear Schr\"odinger
equation (NLS) realized by  the self-induction Vortex  Filament f\/low  (VF). The VF f\/low is a f\/low  in the Euclidean space ${\rm SO}(3)\semip \R^3/{\rm SO}(3)$ (see~\cite{AH}). In \cite{Ha} Hasimoto showed  that, if $u(x,t)$ is a f\/low solution of the VF equation
\begin{gather*}
u_t = \kappa B,
\end{gather*}
where $\kappa$ is the Euclidean curvature of the  curve
$u(\cdot,x)\in \R^3$, $x$ is the arc-length and $B$
is the binormal, then the evolution of  the curvature
and torsion of $u$ is equivalent to the NLS equation via the Hasimoto
transformation $\Phi = \kappa e^{i\int \tau}$. The Hasimoto transformation $(\kappa,\tau)
\to (\nu, \eta)$, with $ \Phi = \nu+i\eta$, is, in fact, induced by
a change from classical Euclidean moving frame  to the  {\it natural
moving frame}  ($\nu$ and $\eta$ are the {\it natural curvatures}, see~\cite{LP1}). Thus, VF is a {\it Euclidean
geometric realization} of NLS if we use natural moving frames. Equivalently, NLS is the {\it invariantization} of~VF. The relation to the Euclidean geometry of the f\/low goes further; consider the evolution
\begin{gather}\label{euinv}
u_t = h T + \frac{h'}{\kappa} N + g B,
\end{gather}
where $\{T,N,B\}$ is the classical Euclidean moving frame and $h$ and $g$
are arbitrary smooth functions of the curvature, torsion and their derivatives. Equation (\ref{euinv})
is the general form of an arc-length preserving  evolution of  space curves,
invariant under the action of the Euclidean group (i.e., $E(n)$ takes solutions to solutions).
Its invariantization can be written as
\begin{gather}\label{keuinv}
\begin{pmatrix} \kappa\\ \tau\end{pmatrix}_t = {\mathcal P} \begin{pmatrix} g\\ h\end{pmatrix},
\end{gather}
where ${\mathcal P}$ def\/ines a Poisson bracket  generated by the {\it second Hamiltonian structure} for NLS
via the Hasimoto transformation, i.e.\
the Hasimoto
transformation is  a {\it Poisson map} (see~\cite{LP1, MSW, MW}). (For more information on inf\/inite dimensional Poisson brackets see \cite{O}, and for more
information on ${\mathcal P}$ see~\cite{MSW}.) Clearly, ${\mathcal P}$ can be generated using
a classical Euclidean moving frame and the invariants $\kappa$ and $\tau$. The NLS equation is a {\it bi-Hamiltonian system}, i.e., Hamiltonian with respect to two compatible Hamiltonian structures. One of the structures
is invertible and a recursion operator can be constructed to generate integrals of the system (see  \cite{Ma} or \cite{O}).
The f\/irst Hamiltonian structure for NLS is invertible and can also be proved to be generated by
the geometry of the f\/low, although it is of a dif\/ferent character as we will see below.
\end{example}

\begin{example}  A second example is that of the Korteweg--de Vries (KdV) equation. If a 1-pa\-rameter family of functions $u(t,\cdot)\in \R$ evolves following the  {\it Schwarzian
KdV equation}
\[
u_t = u' S(u) = u''' - \frac32\frac{(u'')^2}{u'},
\]
where $\displaystyle{S(u) = \frac{u'''}{u'} - \frac32\left(\frac{u''}{u'}\right)^2}$ is the {\it Schwarzian derivative}
of $u$, then $k = S(u)$ itself evolves following the KdV equation
\[
k_t =k''' + 3 kk',
\]
 one of the best known completely
integrable nonlinear PDEs. The KdV equation is Hamiltonian with respect to two
compatible Hamiltonian structures, namely $D=\frac d{dx}$ and $D^3 + 2k D + k'$ (called respectively
f\/irst and second KdV Hamiltonian structures).
As before, if we consider the general curve evolution given by
\begin{gather}\label{up}
u_t = u' h,
\end{gather}
where $h$ is any smooth function depending on $S(u)$ and its derivatives with respect to the parameter
$x$, then $k = S(u)$ evolves following the evolution
\begin{gather}\label{sup}
k_t = (D^3 + 2k D + k') h.
\end{gather}
This time $S(u)$ is the generating {\it differential invariant} associated to the action of ${\rm PSL}(2)$
on~$\RP^1$. That is, any projective dif\/ferential invariant of curves $u(x)$ is a function of $S(u)$ and its derivatives. Equation (\ref{up}) is the most general form for evolutions of reparametrizations of $\RP^1$ (or parametrized ``curves")
  invariant under the action of ${\rm PSL}(2)$.
 Evolution (\ref{sup}) can be viewed as the invariantization of (\ref{up}). Equivalently,
 the family of evolutions in (\ref{up})
 provides   $\R{\rm P}^1$ geometric realizations for the Hamiltonian evolutions
 def\/ined by (\ref{sup}). Thus, one can
obtain geometric reali\-zations in $\RP^1$ not only for KdV, but also for any system which is
Hamiltonian with respect to the second KdV Hamiltonian structure. For example, the Sawada--Koterra
equation
\[
k_t = (D^3 + 2 k D + k') \left(2k'' + \frac12 k^2\right) = 2k^{(5)} + 5 k k''' + 5 k'k'' + \frac52 k'k^2
\]
is bi-Hamiltonian with respect to the same Hamiltonian structures as KdV is. Its Hamiltonian functional
(\ref{sup}) is $h(k) = \int (\frac16 k^3 - (k')^2) dx$. Therefore, the Sawada--Koterra equation has
\[
u_t = u' \left(2S(u)''+\frac12 S(u)^2\right)
\]
as $\RP^1$ realization. (Incidentally, Sawada--Koterra has a second realization as an equi-af\/f\/ine
f\/low, see \cite{O1}.) The manifold $\RP^1$ is an example of a {\it parabolic homogeneous space},
 i.e., a manifold of the form $G/P$ with $G$ semisimple and $P$ a {\it parabolic subgroup}. From
 (\ref{up}) and (\ref{sup}) we can see that the second Hamiltonian structure for KdV can be generated
 with the sole knowledge of $u'$ (a classical projective moving frame along $u$) and $S(u)$, its projective dif\/ferential invariant. The f\/irst KdV structure is also similarly generated, although, again,
 it is of a
 dif\/ferent nature.
 \end{example}

These two simple examples illustrate the close relationship between
the classical geometry of curves and bi-Hamiltonian completely integrable PDEs.
In the last  years many examples of geometric realizations for most
known completely integrable systems have been appearing in the literature. Some
are linked to the geometric invariants of the f\/low
(see for example \cite{A1, A2, DSa, Fe, GR, KQ1, KQ2, LP1,
LP2, M1, M6, MSW, MW, SW, TT, TU1, TU2, YS}).
This list is, by no means, exhaustive as this paper is not meant to be an exhaustive review of the
subject.

 Perhaps the simplest way to understand the close relationship between dif\/ferential invariants and integrable
 systems is through the AKNS representation on one hand and group-based moving frames on the other. If $G$ is a Lie group, a $G$-AKNS representation of a nonlinear PDE
 \begin{gather}\label{PDE}
 k_t = F(k, k_x, k_{xx}, \dots)
 \end{gather}
 is a linear system
 of equations
 \begin{gather*}
 \varphi_x  =  A(t,x,\lambda) \varphi,
\\
 \varphi_t  =  B(x,t,\lambda)\varphi,\qquad
  \varphi(t,x,\lambda)\in G, \qquad A(x,t,\lambda), B(x,t,\lambda) \in \g
\end{gather*}
such that the compatibility condition for the existence of a solution,
 \begin{gather*}
 A_t = B_x + [B,A],
 \end{gather*}
 is independent of $\lambda$ and equivalent to the nonlinear PDE (\ref{PDE}). Such a representation is a basis
 for generating solutions and integrating the system. Indeed, most integrable systems have an AKNS representation. Geometrically, this is thought of as having a $2$-parameter f\/lat connection def\/ined by  $-\frac{d}{dx} + A$ and $-\frac{d}{dt} + B$ along the f\/low (see \cite{HSW}). The bridge to dif\/ferential invariants and dif\/ferential geometry a-la-Cartan appears when one realizes that this $2$-parameter connection is a reduction of the {\it Maurer--Cartan connection} of  $G$ along the f\/low $\varphi$ and the AKNS system could be interpreted as the Serret--Frenet equations and the $t$-evolution of a group-based (right) moving frame  ($\varphi$) along a f\/low
 $u:\R^2\to G/H$ in a certain homogeneous space. This is explained in the  next section.

In this paper we describe how the background geometry of {\it affine} and some {\it  symmetric manifolds} generates Hamiltonian structures and geometric realizations for some completely integrable
systems. Our af\/f\/ine manifolds will be homogeneous manifolds of the form $G\semip \R^n/G$ with $G$ semisimple. They include Euclidean, Minkowski, af\/f\/ine, equi-af\/f\/ine and symplectic geometry among others.
We will also discuss related geometries, like the centro-af\/f\/ine or geometry of star-shaped curves, for which the action of the group is linear instead of af\/f\/ine. On the other hand, our
symmetric manifolds are locally equivalent to a homogeneous manifold of the form $G/H$ where
$\g$, the Lie algebra associated to $G$, has a gradation of the form $\g_{-1}\oplus\g_0\oplus\g_{1}$
and where $g_0\oplus\g_1 = \h$ is the Lie algebra of $H$. These includes projective geometry
($G = {\rm PSL}(n+1)$), the Grassmannian ($G = {\rm SL}(p+q)$), the
conformally f\/lat M\"obius sphere ($G = O(n+1,1)$), the Lagrangian Grassmannian ($G = {\rm Sp}(2n)$), the manifold of reduced pure spinors ($G = O(n,n)$) and more.  We will see how the Cartan geometry of
curves in these manifolds
induces a Hamiltonian structure on the space of dif\/ferential invariants. In the last part of the paper we look closely at the case of symmetric spaces.
We def\/ine {\it differential invariants of projective type} as those gene\-ra\-ted by the action of the group
on second order frames. We then describe how, in most cases, the reduced Hamiltonian structure can be further
reduced or restricted to the space of curves with vanishing
non-projective dif\/ferential invariants. On this manifold the Hamiltonian structure has a (geometrically def\/ined)  compatible Poisson companion. They def\/ine a~bi-Hamiltonian pencil for  some integrable equations of
KdV-type, and they provide geometric realizations for them.
We f\/inally state a conjecture by M. Eastwood on what the presence of these f\/lows might say about the
geometry of curves in symmetric spaces.

  As group based moving frames are relatively new, our next section will describe them in detail. We will also describe their role in AKNS representations.

\section{Moving frames}

The classical concept of moving frame
was developed by \'Elie Cartan \cite{C1, C2}. A classical moving frame along a curve
in a manifold $M$  is a curve in the frame bundle of the manifold over the curve, invariant under
the action of the transformation group under consideration. This method is a very powerful tool, but its
explicit application relied on intuitive choices that were not clear on a general setting.
Some ideas in Cartan's work and later work of Grif\/f\/iths \cite{Gri}, Green \cite{Gre} and others laid the
 foundation for the concept of a group-based moving frame, that is, an equivariant map between
 the jet space of curves in the manifold and the group of transformations. Recent work by
 Fels and Olver \cite{FO1, FO2} f\/inally gave the precise def\/inition of the group-based moving
 frame and extended its application beyond its original geometric picture to an astonishingly
 large group of applications.  In this section we will describe Fels and Olver's moving frame
and its relation to the classical moving frame. We will also introduce
 some def\/initions that are useful to the study of  Poisson brackets  and
bi-Hamiltonian nonlinear PDEs. From now on we will assume $M = G/H$ with $G$
 acting on $M$ via left multiplication on representatives of a class. We will also assume that curves in $M$ are {\it parametrized}
 and, therefore, the group $G$ does not act on the parameter.

 \begin{definition}
Let $J^k(\R,M)$ the space of $k$-jets of  curves, that is,  the set of equivalence classes of curves in
$M$ up to $k^{\rm th}$ order of contact. If we denote by $u(x)$ a curve in $M$ and by
$u_r$ the $r$ derivative of $u$ with respect to the parameter $x$, $u_r = \frac{d^r u}{dx^r}$, the jet space has local
coordinates that can be represented by $u^{(k)}= (x, u, u_1, u_2, \dots, u_k)$.
The group $G$ acts naturally on parametrized curves, therefore it acts naturally
on the jet space via the formula
\[
g\cdot u^{(k)} = (x, g\cdot u, (g\cdot u)_1, (g\cdot u)_2, \dots ),
\]
where by $(g\cdot u)_k$ we mean the formula obtained when one dif\/ferentiates $g\cdot u$ and
then writes the result in terms of
$g$, $u$, $u_1$, etc. This is usually called the {\it prolonged} action of $G$ on $J^k(\R,M)$.
 \end{definition}

\begin{definition}  A function
\[
I: J^k(\R,M) \to \R
\]
is called a $k$th order {\it differential invariant} if it is invariant with respect to the prolonged action of $G$.
\end{definition}

\begin{definition}
A map
\[
\rho: J^k(\R,M) \to G
\]
is called a left (resp. right)  {\it moving frame} if it is equivariant with respect to the prolonged action of $G$
on $J^k(\R,M)$ and the left (resp. right) action of $G$ on itself.
\end{definition}

If a group acts (locally) {\it effectively on subsets}, then for $k$ large enough the prolonged action is locally free on regular jets. This guarantees the existence of a moving frame on a neighborhood of a regular jet (for example, on a neighborhood of a generic curve, see \cite{FO1, FO2}).

 The group-based moving frame already appears in a familiar method for calculating the curvature of a curve $u(s)$ in the
Euclidean plane. In this method one uses a translation to take~$u(s)$ to the origin, and a rotation to make one of the axes tangent to the curve. The curvature can classically be found as the
 coef\/f\/icient of the second order term in the expansion of the curve around $u(s)$. The crucial observation made by Fels and Olver is that {\it the element of the group} carrying out the translation and rotation depends on $u$ and its derivatives and so it def\/ines a map from the jet space to the group. {\it This map is  a right moving frame},
 and it carries all the geometric information of the curve. In fact, Fels and Olver developed
 a similar normalization process to f\/ind right moving frames (see \cite{FO1, FO2} and our next theorem).

\begin{theorem}[\cite{FO1, FO2}]  \label{norm} Let $\cdot$ denote the prolonged
action of the group on $u^{(k)}$ and assume we have {\bfseries normalization equations} of the form
\[
g \cdot u^{(k)} = c_k,
\]
where $c_k$ are constants (they are
called {\bfseries normalization constants}).
Assume we have enough normalization equations so as to determine $g$ as a function of $u, u_1, \dots$. Then $g = \rho$ is a {\bfseries right invariant moving frame}.
\end{theorem}
The direct relation between classical moving frames and group-based moving frames is stated
in the following theorem.

\begin{theorem}[\cite{M2}] \label{classical} Let $\Phi_g: G/H \to G/H$ be defined by multiplication by $g$. That is
$\Phi_g([x]) = [g x]$. Let $\rho$ be a group-based left moving frame with $\rho \cdot o = u$
where $o = [H] \in G/H$. Identify~$d\Phi_\rho(o)$ with an element of $GL(n)$, where $n$ is
the dimension of $M$.

Then, the matrix $d\Phi_\rho(o)$ contains in its columns a classical moving frame.
\end{theorem}
This theorem illustrates how classical moving frames are described only by the
action of the group-based moving frame on f\/irst order frames, while the action on higher
order frames is left out. Accordingly, those invariants determined by the action on higher order
frames will be not be found with the use of a classical moving frame.

Next is the equivalent to the classical Serret--Frenet equations. This concept
if fundamental in our Poisson geometry study.

\begin{definition} Consider $K dx$ to be the horizontal component of the pullback of the left (resp. right)-invariant
Maurer--Cartan form of the group $G$ via a group-based left (resp. right)  moving frame $\rho$. That is
\[
K = \rho^{-1} \rho_x \in \g \qquad ({\rm resp.} \quad K = \rho_x\rho^{-1})
\]
($K$ is the coef\/f\/icient matrix of the f\/irst order dif\/ferential equation satisf\/ied by $\rho$). We call $K$ the {\it left (resp. right) Serret--Frenet equations} for the moving frame $\rho$.
\end{definition}

Notice that, if $\rho$ is a left moving frame, then $\rho^{-1}$ is a right moving frame and their Serret--Frenet equations are the negative of each other. A complete set of
generating dif\/ferential invariants
 can always be found among the
coef\/f\/icients of group-based Serret--Frenet equations, a~crucial dif\/ference with the classical picture. The following
theorem is a direct consequence of the results in \cite{FO1, FO2}. A more general result can be found in
 \cite{Hu}.

\begin{theorem}  Let $\rho$ be a (left or right) moving frame
 along a curve $u$. Then, the coefficients of the (left or right) Serret--Frenet equations
for $\rho$ contain a basis for the space of differential  invariants of the curve. That is, any other differential
invariant for the curve is a function of the entries of $K$ and their derivatives with respect to $x$.
\end{theorem}

\begin{example} Assume $G = \PSL(2)$ so that $M = \RP^1$. The action of $G$ on $\RP^1$
is given by fractional transformations. Assume $\rho = \begin{pmatrix} a&b\\ c&d\end{pmatrix} \in G$ is a (right) moving frame satisfying the normalization equations
\begin{gather*}
\rho \cdot u = \frac{au+b}{cu+d}= 0,\\
\rho \cdot u_1 = \frac{au_1}{cu+d} - \frac{(au+b)cu_1}{(cu+d)^2}= 1, \\
\rho\cdot u_2 = \frac{au_2}{cu+d}- 2 \frac{acu_1^2}{(cu+d)^2} +\frac{au+b}{(cu+d)^3}\left(cu_2(cu+d) + 2 \frac{c^2u^2_1}{(cu+d)^3}\right) = 2\lambda.
\end{gather*}
 Then it is straightforward to check that $\rho$ is completely determined to be
\begin{gather*}
\rho =  \begin{pmatrix}  1& 0\\ \frac12 \frac{u_2}{u_1} - \lambda&1 \end{pmatrix}
\begin{pmatrix} u_1^{-1/2} & 0 \\ 0& u_1^{1/2}\end{pmatrix}\begin{pmatrix} 1 & -u \\ 0&1\end{pmatrix}.
\end{gather*}

A moving frame satisfying this normalization will have the following right Serret--Frenet equation
\begin{gather}\label{KdVSF}
\rho_x =  \begin{pmatrix} -\lambda & -1\\ \frac12 S(u) + \lambda^2 &\lambda\end{pmatrix} \rho.
\end{gather}
This equation is {\it gauge} equivalent to the $\lambda = 0$ equation via the constant {\it gauge}
\begin{gather}\label{gauge}
g = \begin{pmatrix} 1&0\\ \lambda& 1\end{pmatrix}.
\end{gather}
 This gauge $g$ will take the
second normalization constant to zero.

Furthermore, if $u$ is a solution of (\ref{up}), it is known (see \cite{M1}) that the $t$-evolution induced
on $\rho$ is given by
\begin{gather}\label{KdVN}
\rho_t = \begin{pmatrix} -\frac12 h_x - \lambda h& -h\vspace{1mm}\\ \frac12 h_{xx}+\lambda h_x + \lambda^2 h +\frac12 S(u) h& \frac12 h_x+\lambda h\end{pmatrix} \rho.
\end{gather}
\end{example}

We can now see the link between the AKNS representation of KdV and the
evolution of a right moving frame. Assume a completely integrable system (\ref{PDE}) has a geometric realization
which is invariant under the action of the geometric group $G$
\begin{gather}\label{uinv}
u_t = f(\lambda, u, u_1, u_2, \dots).
\end{gather}
Then, under regularity assumptions of the f\/low, the invariantization of (\ref{uinv}) is the integrable system (\ref{PDE}). A right moving frame along $u$ will be a solution of its Serret--Frenet equation
\[
\rho_x = K(t,x,\lambda) \varphi
\]
and the time evolution will induce a time evolution on $\rho$ of the form
\[
\rho_t = N(t,x,\lambda) \varphi.
\]
Furthermore, since (\ref{uinv}) is invariant under the group, both $K$ and $N$ will depend on
the dif\/ferential invariants of the f\/low. These equations are def\/ined by the horizontal
component of the pullback of the  Maurer--Cartan form of the group, $\omega = dg g^{-1}$ by the moving frame
$\rho$, that is $Kdx + Ndt$. If we now evaluate the {\it structure equation} for the Maurer--Cartan form, i.e.,
 $d\omega + \frac12 [\omega,\omega] = 0$,
along $\rho_x$ and $\rho_t$, we get
\[
K_t = N_x + [N,K]
\]
which is exactly the invariantization of the f\/low (\ref{uinv}); therefore it is independent of $\lambda$. Hence,
a~$\lambda$-dependent geometric realization of an integrable system provides  an AKNS representation of the system. See \cite{CIM} for more information.
\begin{example}\label{KdVAKNS} The AKNS representation for KdV is very well known. It is given by the system
\begin{gather*} \varphi_x = \begin{pmatrix} -\lambda & -1 \\ -q & \lambda\end{pmatrix}\varphi,\\
\varphi_t = \begin{pmatrix} -\frac12 q_x - \lambda q + 2\lambda^3 &-q + 2\lambda^2 \\ \frac12 q_{xx}+\lambda q_x +q(-q+2\lambda^2) &\frac12 q_x +\lambda q - 2\lambda^3\end{pmatrix}\varphi.
\end{gather*}

Comparing it to (\ref{KdVSF}) and (\ref{KdVN}) we see that $q = \frac12 S(u)+\lambda^2$ and hence
$u$ will depend on $\lambda$. Furthermore, $h = q-2\lambda^2$ provides $\lambda$-dependent $\RP^1$ geometric realizations for KdV, namely
\[
u_t = u_x\left(-\frac12 S(u) - 3\lambda^3\right).
\]
A complete description of this example can be found in~\cite{CIM}. (Notice that the KdV equation
they represent is dif\/ferent, but equivalent, to our introductory example. This is merely due to
a~dif\/ferent choice of invariant.)
\end{example}

In this paper we will not focus on the study of solutions (see \cite{CIM} instead) but rather on the interaction between geometry and integrable systems. Hence we will largely ignore the spectral parameter $\lambda$ and its role.

\section[Hamiltonian structures generated by group-based moving frames]{Hamiltonian structures generated\\ by group-based moving frames}

Consider the group of loops $\Lo G = C^\infty(S^1, G)$ and its Lie algebra $\Lo \g = C^\infty(S^1,\g)$.
Assume $\g$ is semisimple.
One can def\/ine two natural Poisson brackets on $\Lo\g^\ast$ (see \cite{Se} for more information),
 namely, if $\HH, \F: \Lo\g^\ast \to \R$ are two functionals def\/ined on $\Lo\g^\ast$ and if
 $L \in \Lo\g^\ast$, we def\/ine
 \begin{gather}\label{br1}
 \{\HH, \F \}_1(L) = \int_{S^1} \left\langle \left(\var{\HH} \right)_x + {\rm ad}^\ast\left(\var{\HH}\right) (L) , \var \F\right\rangle dx,
 \end{gather}
where $\langle \ , \ \rangle$ is the natural coupling between $\g^\ast$ and~$\g$, and where
$\var{\HH}$ is the variational derivative of $\HH$ at $L$ identif\/ied, as usual, with an element of
$\Lo\g$.

One also has a compatible family of second brackets, namely
\begin{gather}\label{br2}
\{\HH, \F\}_2(L) = \int_{S^1} \left\langle {\rm ad}^\ast\left(\var{\HH}\right) (L_0) , \var \F\right\rangle dx,
\end{gather}
where $L_0\in\g^\ast$ is any constant element.
Since $\g$ is semisimple we can identify $\g$ with its dual~$\g^\ast$ and we will do so from now on.

From now on we will also assume that our curves on homogeneous manifolds have a {\it group monodromy}, i.e., there exists $m\in G$ such that
\[
u(t+T) = m\cdot u(t),
\]
where $T$ is the period. Under these assumptions, the Serret--Frenet equations will be periodic. One could, instead, assume that~$u$ is asymptotic at~$\pm \infty$, so that the invariants will vanish at
inf\/inity. We would then work with the analogous of~(\ref{br1}) and~(\ref{br2}).

The question we would like to investigate next is whether or not these two brackets can be reduced
to the space of dif\/ferential invariants, or the space of dif\/ferential invariants associated to special types of f\/lows. We will describe af\/f\/ine and symmetric cases separately.

\subsection[Affine manifolds]{Af\/f\/ine manifolds}

Assume $M = (G\semip\R^n)/G$ is an af\/f\/ine manifold, $G$ semisimple. In this case a moving frame
can be represented as
\begin{gather}\label{affrho}
\rho = \begin{pmatrix} 1 & 0 \\ \rho_u & \rho_G\end{pmatrix}
\end{gather}
acting on $\R^n$ as
\[
\rho \cdot u = \rho_G u + \rho_u.
\]
A left invariant moving frame with $\rho\cdot o=  u$ will hold $\rho_u = u$ and, in view of  Theorem
\ref{classical}, $\rho_G$ will have in
its columns a classical moving frame. In this case $K = \rho^{-1}\rho_x$ is given by
\[
K = \begin{pmatrix} 0&0\\ \rho^{-1}_G (\rho_u)_x & \rho^{-1}_G (\rho_G)_x \end{pmatrix}.
\]
In \cite{M2} it was shown that $\rho_G^{-1} (\rho_u)_x$ contains all f\/irst order dif\/ferential
invariants. It was also explained how one could make this term constant by choosing
a special parametrization if necessary. Let's call that constant $\rho_G^{-1} (\rho_u)_x = \Lambda$. Our main tool
to f\/ind Poisson brackets is via reduction, and as a previous step, we need to write the
space of dif\/ferential invariants as a~quotient in $\Lo\g^\ast$. The proof of the following
Theorem  can be found in \cite{M2}.
\begin{theorem}[\cite{M2}] Let $N\subset G$ be the isotropy subgroup of $\Lambda$. Assume that we choose
moving frames as above and let $\K$ be the space of Serret--Frenet equations determined by these
moving frames for curves in a neighborhood of a generic curve $u$. Then, there exists an open set of $\Lo\g^\ast$, let's call it $U$, such that $U/\Lo N \cong \K$, where $\Lo N$ acts on $\Lo\g^\ast$ using the
gauge (or Kac--Moody) transformation
\begin{gather}\label{km}
a^\ast(n)(L) = n^{-1}n_x + n^{-1} L n.
\end{gather}
\end{theorem}

In view of this theorem, our next theorem comes as no surprise.
\begin{theorem}[\cite{M2}] The Hamiltonian structure \eqref{br1} reduces to $U/\Lo N \cong \K$ to define
a Poisson bracket in the space of differential invariants of curves.
\end{theorem}

\begin{example} If we choose $G = SO(3)$ and $M$ the Euclidean space,  for appropriate
choice of normalization constants our left moving frame is given by
\[
\rho = \begin{pmatrix} 1&0\\ u& T~ N~ B\end{pmatrix},
\]
where $\{T,N,B\}$ is the classical Euclidean Serret--Frenet frame. Its Serret--Frenet
equations will look like
\[
K = \rho^{-1}\rho_x = \begin{pmatrix} 0&0&0&0\\ (u_1\cdot u_1)  & 0&-\kappa&0\\ 0&\kappa & 0 &-\tau\\ 0&0&\tau&0\end{pmatrix}.
\]
In this case, if we choose to parametrize our curve by arc-length, $\Lambda = e_1$ where, as usual, we denote by $e_k$ the standard basis of $\R^n$. The matrix  $K$ can be clearly identif\/ied with its $\oa(3)$ block and hence $\K$ can be considered as a subspace of $\Lo\oa(3)^\ast$. The isotropy subgroup~$N$ is given by matrices of the form $\begin{pmatrix} 1 & 0\\ 0 & \Theta\end{pmatrix}$ with $\Theta \in SO(2)$.

Using this information we can f\/ind the reduced bracket algebraically. For this we take any functional
$h:\K \to \R$. Let's call $\HH$ an extension of $h$ to $\Lo\oa(3)^\ast$, constant on the gauge leaves
of~$\Lo N$. Its variational derivative at $K$ needs to look like
\[
\frac{\delta \HH}{\delta L}(K) = \begin{pmatrix} 0&\hk&\alpha\vspace{1mm}\\ -\hk&0&\hta\vspace{1mm}\\ -\alpha & -\hta&0\end{pmatrix}
\]
for some $\alpha$ to be determined. Since $\HH$ is constant on the gauge leaves of $\Lo N$
\[
\left\langle n^{-1}n_x + n^{-1} K n, \frac{\delta \HH}{\delta L}(K)\right\rangle = 0,
\]
for any $n \in \Lo N$. This is equivalent to
\[
\left(\frac{\delta \HH}{\delta L}(K)\right)_x + \left[ K, \frac{\delta \HH}{\delta L}(K)\right] \in \Lo \n^o,
\]
where $\n$ is the Lie algebra of $N$ and $\n^o$ is its annihilator. From here
\[
\begin{pmatrix} 0&(\hk)_x - \alpha \tau & \alpha_x - \kappa \hta+\tau\hk\vspace{1mm}\\ \ast & 0 & (\hta)_x + \kappa \alpha\vspace{1mm}\\ \ast&\ast&0\end{pmatrix} = \begin{pmatrix} 0&\ast&\ast\\ \ast&0&0\\ \ast&0&0\end{pmatrix},
\]
where $\ast$ indicates entries that are not, at least for now, relevant. Hence $\alpha = -\frac 1\kappa(\hta)_x$. The bracket is thus given by
\begin{gather*}
\{f,h\}^R_1(K)  =
\int_{S^1}{\rm tr} \left(\begin{pmatrix} 0&\ast&\ast\vspace{1mm}\\ -(\hk)_x-\frac\tau\kappa(\hta)_x& 0 & 0\vspace{1mm}\\  \left(\frac1\kappa\left(\hta\right)_x\right)_x + \kappa \hta-\tau\hk &0&0\end{pmatrix} \begin{pmatrix} 0&\fk&-\frac1\kappa(\ft)_x\vspace{1mm}\\ -\fk&0&\ft\vspace{1mm}\\ \frac1\kappa(\ft)_x&-\ft&0\end{pmatrix}\right) dx
\\
\phantom{\{f,h\}^R_1(K)  }{} =  -2 \int_{S^1} \begin{pmatrix} \fk &\ft \end{pmatrix} {\mathcal R} \begin{pmatrix} \hk \vspace{1mm}\\ \hta\end{pmatrix} dx,
\end{gather*}
where ${\mathcal R}$ is
\[
{\mathcal R} = \begin{pmatrix} D& \frac\tau\kappa D\vspace{1mm}\\ D\frac\tau\kappa& -D-D\frac1\kappa D\frac1\kappa D\end{pmatrix}.
\]
The second Hamiltonian structure (\ref{br2}) can also be reduced to $\K$ with the general  choice $L_0 = \begin{pmatrix} 0& a&b\\ -a&0&c\\ -b&-c&0\end{pmatrix}$. The reduced bracket is found when applying
(\ref{br2}) to the variational derivatives of  extensions that, as before, are constant on the $\Lo N$ leaves. Thus, it is straightforward to check that
the second reduced bracket is given by
\[
\{f,h\}^R_2(K) = 2 \int_{S^1} \begin{pmatrix} \fk&\ft\end{pmatrix}\left(a {\mathcal A} + b{\mathcal B} + c {\mathcal C}\right) \begin{pmatrix} \hk\vspace{1mm}\\ \hta\end{pmatrix} dx,
\]
where
\[
{\mathcal A} = \begin{pmatrix} 0&0\vspace{1mm}\\ 0& \frac1\kappa D - D\frac1\kappa\end{pmatrix},\qquad {\mathcal B} = \begin{pmatrix} 0&1\\ -1&0\end{pmatrix},\qquad {\mathcal C} = \begin{pmatrix} 0&\frac1\kappa D\vspace{1mm}\\ D\frac1\kappa&0\end{pmatrix}.
\]
These are all Hamiltonian structures and they appeared in \cite{MSW}. In fact, the structure ${\mathcal P}$ shown in the introduction can be written as ${\mathcal P} = -{\mathcal R} {\mathcal C}^{-1} {\mathcal R}$ and, hence, ${\mathcal P}$ is in the Hamiltonian hierarchy generated by ${\mathcal R}$ and ${\mathcal C}$. A study of integrable systems associated to these brackets, and their geometric realizations, was done in \cite{MSW}. See also \cite{I}.
\end{example}

Our f\/irst reduced bracket is directly related to geometric realizations. In fact, by choosing a~special parameter $x$ we can obtain geometric realizations of systems that are Hamiltonian with respect to the reduced bracket. We explain this next. Let~$\rho$ be given as in (\ref{affrho}) and assume~$u$ is a
solution of the invariant equation
\begin{gather}\label{uev}
u_t = \rho_G {\bf r},
\end{gather}
where ${\bf r} = (r_i)$ is a dif\/ferential invariant vector, that is, $r_i$ are all functions of the entries
of $K$ and its derivatives. If the parameter has been f\/ixed so as to guarantee that $\rho_G^{-1}\rho_u = \Lambda$ is constant, then ${\bf r}$ has to be modif\/ied to guarantee that the evolutions (\ref{uev}) preserve the parameter. (In the running example $\rho_G = (T, N, B)$ and $\rho_G{\bf r} = r_1T+r_2N+r_3B$ with $r_2 = \frac{r_1'}{\kappa}$ once the arc-length is chosen as parameter.)

\begin{theorem}[\cite{M2}] If there exists a  Hamiltonian $h:\K \to \R$ and a local extension $\HH$
constant on the leaves of $N$ such that
\begin{gather}\label{grc}
\frac{\delta\HH}{\delta L}(K) \Lambda = \rr_x + K\rr,
\end{gather}
then the invariantization  of evolution \eqref{uev} is Hamiltonian with respect to
the reduced bracket $\{\ , \ \}^R_1$ and its associated Hamiltonian is $h$.
\end{theorem}

If (after choosing a special parameter if necessary)  relation (\ref{grc}) can be solved for $\rr$
given a certain Hamiltonian $h$, then the Theorem guarantees a geometric realization for the reduced
Hamiltonian system. Such is the case for the VF f\/low, Sawada--Koterra \cite{M2, O1},
modif\/ied KdV and others~\cite{MSW}.

One can f\/ind many geometric realizations of integrable systems in af\/f\/ine manifolds (see, for example, \cite{A2, DSa, I, KQ1,KQ2, LP1, SW, TT, TU2, YS}). Many of these are realizations of modif\/ied KdV equations (or its generalizations), sine-Gordon and Schr\"odinger
f\/lows. These systems have geometric realizations also in non-af\/f\/ine manifolds (see \cite{A1, KQ1,KQ2, TT, TU1,TU2}). Nevertheless, a~common feature to the generation of these realizations is the existence of a classical moving frame that resembles the classical {\it natural} moving frame, that is, the derivatives of the non-tangential vectors of the classical frame all have a tangential direction. Thus, it seem
to be the case that the existence  of geometric realizations for these systems is linked to the existence of a natural frame. This close relationship between geometry and the type of integrable
system is perhaps clearer in our next study, that of symmetric manifolds.

Before moving on, we have one f\/inal comment in this line of thought.
 There are other manifolds whose geometry is given by a linear (rather than af\/f\/ine) action of the group.
 We can still follow a similar approach, reduce the brackets and study Hamiltonian structures on the
 space of dif\/ferential invariants. For example, in the case of centro-af\/f\/ine geometry one considers the linear action of $\SL(n)$ on $\R^n$ and the associated geometry is that of star-shaped curves. If we assume that curves are parametrized by the centro-af\/f\/ine arc-length, we can reduce both brackets and
 obtain a pencil of Poisson brackets. This pencil coincides with the bi-Hamiltonian structure of KdV. Indeed, a geometric realization for KdV was found by Pinkall in \cite{Pi}. This realization is the one
 guaranteed by the reduction, as explained in \cite{CIM}. The interesting aspect of the centro-af\/f\/ine case is the following: there is a natural identif\/ication
 of a~star-shaped curve with a projective curve.  The identif\/ication is given by the intersection of the curve
 with the lines going through the origin. If the star-shaped curve is {\it nondegenerate} (that is, $\det(\gamma, \gamma_x,\dots,\gamma^{(n-1)}) \ne 0$. For example, in the planar case the curve is never in
 the radial direction), the identif\/ication is well-def\/ined. Furthermore, if we parametrized star-shaped
 curves with centroaf\/f\/ine arc-length (that is, if $\det(\gamma, \gamma_x,\dots,\gamma^{(n-1)}) = 1$),
 the identif\/ication is 1-to-1 and the geometries are Poisson-equivalent, the Poisson isomorphism given by the identif\/ication. In fact,
 Pinkall's geometric realization is the star-shaped version of the Schwarzian KdV under this relation
 (see \cite{CIM} for more details).
The existence of a geometric realization for KdV seems to imply the existence of a background projective geometry.

\subsection{Symmetric manifolds}

Assume that $M$ is a  symmetric manifold which is locally equivalent to $G/H$ with: (a)
$G$  semisimple; (b)
its (Cartan) connection  given
by the Maurer--Cartan form (i.e.\  the manifold is f\/lat);
(c)~the Lie algebra $\g$ has a gradation of length 1, i.e.
\begin{gather}\label{grad}
\g = \g_{-1} \oplus \g_0\oplus \g_1
\end{gather}
with $\h = \g_0\oplus\g_1$, where $\h$ is the Lie algebra of $H$.

If $M$ is a symmetric manifold of this type,  $G$ splits locally as
$G_{-1}\cdot G_0\cdot G_1$ with $H$ given by $G_0\cdot G_1$. The subgroup $G_0$
is called {\it the isotropic subgroup} of $G$ and it is the component of $G$ that acts linearly
on $G/H$ (for more information see \cite{Ba} or~\cite{Oc}). That means $G_0$ is the component
of the group acting on f\/irst
order frames. According to Theorem \ref{classical}, the $\rho_0$ factor of a left moving frame $\rho$
will determine a classical moving frame (see also \cite{M1}).

As in the previous case, the basis for the def\/inition of a  Poisson bracket on the space of invariants
is to express that space as a quotient in $\Lo\g^\ast$. This is the result in the following Theorem.
 For a complete description and proofs
see \cite{M1}. Notice that, if $\rho$ is a (right) moving frame along a curve in a
 symmetric manifold with $\rho\cdot u = o$,
then $\rho\cdot u_1$ is always constant. In general $\rho\cdot u_1$ is described by f\/irst order invariants, but
curves in symmetric manifolds do not have non-constant f\/irst order dif\/ferential invariants (invariants
are third order or higher), and hence $\rho\cdot u_1$ must be constant.

\begin{theorem}[\cite{M1}]\label{main}  Let $M = G/H$ be a  symmetric manifold as above.
Assume that for every curve in a neighborhood of a generic curve $u$ in $M$ we choose a left moving frame $\rho$ with $\rho\cdot o = u$
and $\rho^{-1}\cdot u_1 = \hat \Lambda$ constant. Assume that we choose a section of $G/H$ so we can locally identify
the manifold with $G_{-1}$ and its tangent with $\g_{-1}$. Let $\Lambda\in\g_{-1}$ represent $\hat\Lambda$ and let $\K$ be the manifold of  Serret--Frenet
equations for $\rho$ along curves in a neighborhood of $u$.  Clearly $\K \subset \Lo\g^\ast$. Denote by $\langle \Lambda\rangle$ the linear
subspace of $C^\infty(S^1,\g^\ast)$ given by $\langle \Lambda \rangle = \{ \alpha \Lambda, \alpha(x) >0\}$.

 Then the space $\K$   can be described as a quotient $U/ \N$, where $U$ is an open set of
$\langle \Lambda\rangle \oplus \Lo\g_0 \oplus \Lo\g_1$ and where $\N = \N_0\cdot \Lo G_1 \subset \Lo G_0\cdot \Lo G_1$ acts on $U$ via the Kac--Moody action \eqref{km}.
The subgroup~$\N_0$ is the isotropy subgroup of $\langle\Lambda\rangle$ in $\Lo G_0$.
\end{theorem}

As before, after writing $\K$ as a quotient, one gets a reduction theorem.

\begin{theorem}[\cite{M1}] The Poisson bracket \eqref{br1} can be reduced to $\K$ and there exists a
well-defined Poisson bracket $\{ \ ,\ \}^R_1$ defined on a generating set of independent differential invariants.
\end{theorem}

\begin{example}As we saw before, the (left) Serret--Frenet equations for the $\RP^1$ case are given by
\[
K = \begin{pmatrix} 0&1\\ k& 0\end{pmatrix}
\]
where $k = -\frac12 S(u)$. The splitting of the Lie algebra $\sla(2)$ is given by
\[
\begin{pmatrix} 0& \beta\\ 0&0\end{pmatrix} + \begin{pmatrix} \alpha & 0 \\ 0 & -\alpha\end{pmatrix} +
\begin{pmatrix} 0& 0\\ \gamma&0\end{pmatrix} \in \g_{-1} + \g_0+\g_1.
\]
 In this case $\Lambda = \begin{pmatrix} 0&1\\ 0&0\end{pmatrix}\in \g_{-1}$. The isotropic subgroup of $\langle\Lambda\rangle$ in $G_0$ is $G_0$ itself, and so $\N = \Lo G_0\cdot \Lo G_1$ or subspace
 of lower triangular matrices. To
reduce the bracket (\ref{br1}) we need to have a functional $h:\K \to \R$ and to f\/ind an extension $\HH:\M \to \R$
such that
\begin{gather}\label{condi}
\left(\frac{\delta\HH}{\delta L}(K)\right)_x + \left[ K, \frac{\delta\HH}{\delta L}(K)\right] \in \n^0
\end{gather}
for any $K\in \K$, where $\n = \Lo\g_0\oplus \Lo\g_1$ and $\n^0$ is its annihilator (which we can
identify with $\Lo\g_1 = \Lo\g_{-1}^\ast$).  Also, if $\HH$ is an extension of $h$, its variational derivative at $\K$ will be given by
\[
\frac{\delta\HH}{\delta L}(K) = \begin{pmatrix} a& \frac{\delta h}{\delta k}(k) \\ b&-a\end{pmatrix},
\]
where $a$ and $b$ are to be determined. On the other hand condition (\ref{condi}) translates into
\begin{gather*}
\begin{pmatrix} a_x +b-k\frac{\delta h}{\delta k}(k) & \left(\frac{\delta h}{\delta k}(k)\right)_x -2a\vspace{1mm}\\ b_x + 2ka
& -a_x -b+k\frac{\delta h}{\delta k}(k)\end{pmatrix} = \begin{pmatrix} 0&0\\ \ast&0\end{pmatrix}.
\end{gather*}

From here $a = \frac 12 \left(\frac{\delta h}{\delta k}(k)\right)_x$ and $b = k\frac{\delta h}{\delta k}(k) - \frac12 \left(\frac{\delta h}{\delta k}(k)\right)_{xx}$. We are now ready to f\/ind the reduced bracket. If $f,h: \K\to \R$ are two functionals and $\F$ and $\HH$ are extensions vanishing on the $\N$-leaves, the reduced bracket is given by
\begin{gather*}
\{ f, h\}^R_1(K) = \int_{S^1} {\rm tr}\left(\frac{\delta\F}{\delta L}(K)\left\{ \left(\frac{\delta\HH}{\delta L}(K)\right)_x + \left[K, \frac{\delta\HH}{\delta L}(K)\right]\right\}\right) dx
\\
\phantom{\{ f, h\}^R_1(K)}{} = \int_{S^1} {\rm tr} \left(\begin{pmatrix} \ast   & \frac{\delta f}{\delta k}(k)\\ \ast&\ast\end{pmatrix}\begin{pmatrix} 0&0\vspace{1mm}\\ k \left(\frac{\delta h}{\delta k}(k)\right)_x +( k\frac{\delta h}{\delta k}(k))_x - \frac12 \left(\frac{\delta h}{\delta k}(k)\right)_{xxx} & 0\end{pmatrix}\right) dx
\\
\phantom{\{ f, h\}^R_1(K)}{} = \int_{S^1} \frac{\delta f}{\delta k}(k)\left( -\frac12 D^3 + k D + D k\right) \frac{\delta h}{\delta k}(k) dx.
\end{gather*}
The dif\/ference in coef\/f\/icients as compared to the introductory example is due to the fact that $k = -\frac12 S(u)$ and not $S(u)$.

As it happens, the companion bracket also reduces for $L_0 = \Lambda^\ast = \begin{pmatrix} 0&0\\1&0\end{pmatrix}$. Indeed, it is given~by
\begin{gather*}
\{f, h\}_2^R(K)  =  \int_{S^1} {\rm tr}\left\{ \frac{\delta\F}{\delta L}(K) \left[ \begin{pmatrix} 0&0\\ 1&0\end{pmatrix}, \frac{\delta\HH}{\delta L}(K)\right]\right\}
\\
\phantom{\{f, h\}_2^R(K)}{}  =  \int_{S^1} {\rm tr}\left\{ \begin{pmatrix} \frac12 \left( \frac{\delta f}{\delta k}(k)\right)_x& \frac{\delta f}{\delta k}(k)\vspace{1mm}\\ \ast & -\frac12 \left( \frac{\delta f}{\delta k}(k)\right)_x\end{pmatrix}\begin{pmatrix}  -\frac{\delta h}{\delta k}(k) & 0 \vspace{1mm}\\  \left( \frac{\delta h}{\delta k}(k)\right)_x&  \frac{\delta h}{\delta k}(k)\end{pmatrix}\right\}dx
\\
\phantom{\{f, h\}_2^R(K)}{} =  2
\int_{S^1} \frac{\delta f}{\delta k}(k) D \frac{\delta h}{\delta k}(k) dx.
\end{gather*}
\end{example}

It is not true in general that (\ref{br2}) is also reducible to $\K$. In fact, one f\/inds that for $M = \R{\rm P}^n$ and $G = {\rm PSL}(n+1)$ the second bracket (\ref{br2}) is also reducible to $\K$ when $L_0 = \Lambda^{\ast}
\in \g^\ast$. The resulting two brackets are the f\/irst and second Hamiltonian structure for
Adler--Gel'fand--Dikii f\/lows. But if $M$ is the so-called Lagrangian Grassmannian, $G = {\rm Sp}(4)$, the second
bracket is never reducible to $\K$~\cite{M4}.

One interesting comment on the connection to AKNS representations: as before, the reduction of the bracket (\ref{br1}) is  directly linked to geometric realizations. But the reduction of~(\ref{br2}) indicates the existence of an AKNS representation and
an integrable system. In the KdV example we described how
the Serret--Frenet equa\-tions~(\ref{KdVSF}) were gauge equivalent to $\lambda = 0$ using the gauge~(\ref{gauge}). If we gauge the $x$-evolution of the KdV  AKNS representation in Example
\ref{KdVAKNS} by that same element we get that the matrix $A$ changes into
\[
A_\lambda  = \begin{pmatrix} 0 & 1 \\ q-\lambda^2&0\end{pmatrix} = A -\lambda^2 L_0.
\]
 Therefore, up to a constant gauge, $L_0$  indicates the position of the spectral parameter in the KdV AKNS representation in Example~\ref{KdVAKNS}.  In fact, it goes further.  One can prove that the coef\/f\/icient~$h$ (as in~(\ref{up})) of the $\lambda$-dependent realizations for KdV determined by this AKNS representation (that is, $h = q-2\lambda^2$) is given by the variational derivative of the Hamiltonian functional  used to write KdV as Hamiltonian system in the pencil $\{\ ,\ \}^R_1 - \lambda^2 \{\ ,\ \}^R_2$. In the same fashion, the NLS,
 when written in terms of $\kappa$ and $\tau$ as in (\ref{keuinv}), has a second Hamiltonian structure obtained when reducing (\ref{br2}) with the choice $L_0 = \begin{pmatrix} 0&0&0\\ 0&0&1\\ 0&-1&0\end{pmatrix}$. One can see that this system has an AKNS representation with $x$ evolution given by
 \[
 \rho_x = \begin{pmatrix} 0&\kappa&0\\ -\kappa&0&\tau - \lambda \\ 0&-\tau+\lambda & 0\end{pmatrix}\rho
 \]
 and where the $t$ component is determined by the evolution induced on the {\it right} moving frame~$\rho$
 by a $u$ evolution whose invariant coef\/f\/icients $h$ and $g$ as jn (\ref{euinv}) are given by the variational derivative of the Hamiltonian used to write this Euclidean representation of NLS as Hamiltonian with
 respect to the pencil $\{\ ,\ \}^R_1 - \lambda \{\ ,\ \}_2^R$. This must be a known fact on integrable systems, but we could not f\/ind it in the literature. For a complete description of this relation see~\cite{CIM}.

Back to symmetric spaces. Recall that $\g = \g_{-1}\oplus\g_0\oplus\g_1$ with $\h = \g_0\oplus\g_1$, so that we can identify
$T_xM$ with $\g_{-1}$. Recall also that $d\Phi_\rho(o) = (T_1 \dots T_n)$
is a classical moving frame.
\begin{theorem}[\cite{M1}] Assume $\rho$ is as in the statement of Theorem~{\rm \ref{main}}. Assume $u(t,x)$
is a solution of the evolution
\begin{gather}\label{uev*}
u_t = d\Phi_\rho(o) {\bf r} = r_1 T_1 + \dots + r_n T_n,
\end{gather}
where ${\bf r} = (r_i)$ is a vector of differential invariants. Then, if there exists a Hamiltonian functional
$h : \K \to \R$ and an extension of $h$, $\HH: U \to \R$, constant of the leaves of $\N$ and such that $\left[\var{\HH}\right]_{-1} = {\bf r}$ (the subindex $-1$ indicates the component in $\g_{-1}$), then the evolution induced on the generating system of differential
invariants defined by $\K$ is Hamiltonian with associated Hamiltonian $h$.
\end{theorem}

Let us call the generating dif\/ferential invariants ${\bf k}$. Given a Hamiltonian system ${\bf k_t} = \xi_{h}({\bf k}) = \Pm({\bf k}) \frac{\delta h}{\delta {\bf k}}$,  its associated geometric evolution
(\ref{uev*}) will be its geometric realization in $M$. Therefore they always exist, the previous theorem guarantees their
existence. Indeed, one only needs to extend the functional $h$ preserving
the leaves. As explained in \cite{M1} and described in our running examples, one can f\/ind $\frac{\delta \HH}{\delta L}(K)$ along $\K$
explicitly
using a simple algebraic process. Then, the coef\/f\/icients $r_i$ of the realization are given by $\left(\frac{\delta \HH}{\delta L}(K)\right)_{-1}$, as identif\/ied with the tangent to the manifold using our section
of $G/H$. Some examples are given in \cite{M1, M3, M4, M5}.

As we have previously pointed out,  the second bracket does never reduce to the space of
invariants when $G = \Sp(4)$, the Lagrangian Grassmannian. Still, one can select a submanifold
of invariants and study reductions of the brackets on the submanifold where the other invariants
vanish. This is equivalent to studying Hamiltonian evolutions of special types of f\/lows on $\Sp(4)/H$.
In particular, the author def\/ined dif\/ferential invariants of projective type in \cite{M3}. She then
showed how both brackets (\ref{br1}) and (\ref{br2}) (for some choice of $L_0$) could be reduced on f\/lows of curves with
vanishing
non-projective dif\/ferential invariants. The reductions produced Hamiltonian structures and geometric
realizations for integrable systems of KdV-type. This, and its implication for the geometry of $G/H$
is described next.

\section[Completely integrable systems of KdV type associated to differential invariants of projective type]{Completely integrable systems of KdV type associated\\ to dif\/ferential invariants of projective type}

There are some dif\/ferential invariants of curves in  symmetric spaces that one might call
{\it of projective type}. They are generated by the action of the group on second frames (hence
they cannot be found using a classical moving frame) and most of them closely resemble the
Schwarzian derivative.  In terms of the gradation
$\g = \g_{-1}\oplus \g_0\oplus\g_1$, if we choose an appropriate moving frame, these invariants will
 appear in $K_1$, where $K = \rho^{-1}\rho_x = K_{-1}+ K_0+ K_1$ is the graded splitting of the Serret--Frenet equation associated
to the moving frame $\rho$.

Of course, the simplest example of dif\/ferential invariants of projective-type are projective dif\/ferential invariants. As it was shown in \cite{M1}, a moving frame can be chosen so that all dif\/ferential invariants
appear in the $\g_1$ component of the Serret--Frenet equations, while all entries outside $\g_1$
are constant. More examples of dif\/ferential invariants of projective type appear in \cite{M3, M4, M5} and \cite{M7}.

In this series of papers the author showed how in many symmetric  spaces  one can f\/ind
geometric realizations inducing evolutions of KdV type on the dif\/ferential invariants of
projective type. Indeed, if $G/H$ is a symmetric
space as above, it is known \cite{K} that $\g$ is the direct sum of the following simple Lie
algebras:
\begin{enumerate}\itemsep=0pt
\item\label{proj} $\g = \sla(p+q)$ with $p,q\in\Z^+$. If $q = 1$ then $G/H \equiv \RP^n$. In general $G/H$ is the {\it Grassmannian}.
 \item\label{Lag} $\g = \spa(2n)$, the manifold $G/H$ is called the {\it Lagrangian Grassmannian} and it can be
 identif\/ied with the manifold of Lagrangian planes in $\R^{2n}$.
 \item\label{spi} $\g= \oa(n,n)$, the manifold $G/H$ is called the manifold of {\it reduced pure spinors}.
\item\label{conf} $\g = \oa(p+1,q+1)$ with $p,q\in \Z^+$. If $q = 0$ the manifold $G/H$ is isomorphic to the
 {\it M\"obius sphere}, the local model for f\/lat conformal manifolds.
 \item\label{exc} Two exceptional cases, $\g = E_6$ and $\g = E_7$.
\end{enumerate}
We will next describe the situation for each one of the \ref{proj}--\ref{conf} cases above. The Grassmannian case~(\ref{proj}) (other than $q = 1$) and the exceptional cases (\ref{exc}) have not yet been
studied.

\subsection{Projective case}

Let $G = {\rm PSL}(n+1)$. If $g \in G$ then, locally
\[
g  = g_{-1}g_0g_1= \left(\begin{array}{cc}I&u\\ 0&1\end{array}\right) \left(\begin{array}{cc}\Theta&0\\ 0&(\mathrm{det}\Theta)^{-1}\end{array}\right) \left(\begin{array}{cc}I&0\\ v^T&1\end{array}\right)
\]
with $u,v\in \R^n$ and $\Theta\in {\rm GL}(n)$. If we def\/ine $H$ by the choice $u = 0$, then $G/H\cong \RP^n$, the $n$-projective space. This factorization corresponds to the splitting given by the gradation (\ref{grad}). A section for $G/H$ can be taken to be the $g_{-1}$ factor. As with any other homogeneous space, the action of $G$ on this section is completely determined by the relation
\[
g \begin{pmatrix} I & u\\ 0&1\end{pmatrix} = \begin{pmatrix} I & g \cdot u\\ 0&1\end{pmatrix} h
\]
for some $h\in H$.

The corresponding splitting of the Lie algebra is given by
\[
V = V_{-1}+V_0+V_1 = \begin{pmatrix} 0&a\\ 0&0\end{pmatrix} +   \begin{pmatrix} A&0\\ 0&-{\rm tr}A\end{pmatrix} + \begin{pmatrix} 0&0\\ b^T&0\end{pmatrix} \in \g_{-1}\oplus\g_0\oplus\g_1.
\]
We can identify the f\/irst term above with the tangent to the manifold.

The following theorem describes the type of moving frame and invariant manifold we will choose.
It is essential to choose a simple enough representation for $\K$ so that one can readily recognize
the result of the reduction. Naturally, any choice of moving frame will produce a~choice for $\K$ and
a Hamiltonian structure. But showing the equivalence of Poisson structures is a non-trivial problem,
and hence its recognition is an important part of the problem. As we said before, all dif\/ferential
invariants are of projective-type.

 \begin{theorem}[reformulation of Wilczynski \cite{W}] There exists a
left moving frame $\rho$ along nondegenerate curves in $\R\mathrm{P}^n$ such that its Serret--Frenet equations are defined by matrices of the form
\[
K =  \left(\begin{array}{ccccc} 0 & 1&0&\dots &0\\0&0&1&\dots &0\\
\vdots&\vdots&\ddots&\ddots&\vdots\\ 0&0&\dots&0&1\\ k_1&k_2&
\dots&k_n&0\end{array}\right),
\]
where $k_i$, $i=1,\dots,n$ are, in general, a generating combination of the \emph{Wilczynski projective invariants} and their derivatives.
\end{theorem}

For the precise relation between these invariants and Wylczynski's invariants,  see \cite{CIM}.
The following result was originally obtained by Drinfel'd and Sokolov in \cite{DS}. Their description of the quotient is not the
same as ours, but \cite{M1} showed that our reduction and theirs are equivalent.
\begin{theorem}[\cite{DS}] Assume $\K$ is represented by matrices of the form above. Then, the reduction of \eqref{br1} to $\K$ is given by the
Adler--Gel'fand--Dikii  (AGD) bracket or second Hamiltonian structure for generalized KdV. The bracket \eqref{br2}
reduces for the choice $L_0 = e_n^\ast$ and its reduction is the first Hamiltonian structure for generalized KdV
equations.
\end{theorem}

Finally,  the following theorem identif\/ies geometric realizations for any f\/low Hamiltonian with respect to
the reduction of (\ref{br1}); in particular, it provides geometric realizations for generalized KdV equations
or AGD f\/lows. The case $n=2$ was originally proved in \cite{HLM}.
\begin{theorem}[\cite{M6}] Assume $u:J\subset \R^2 \to {\rm PSL}(n+1)/H$ is a solution
 of
 \[
 u_t = h_1 T_1 + h_2 T_2 + \dots + h_nT_n,
 \]
 where $T_i$ form a {\it projective classical moving frame}.

 Then, ${\bf k} = (k_i)$ satisfies an equation of the form
  \[
 {\bf k}_t = P {\bf h},
 \]
 where ${\bf k} = (k_1,\dots,k_n)^T$,  ${\bf h} = (h_1,\dots,h_n)^T$ and where $P$ is the Poisson
 tensor defining  the {\bfseries Adler--Gel'fand--Dikii Hamiltonian structure}.
 In particular, we obtain a projective geometric realization for
a {\bfseries generalized KdV system} of equations.
  \end{theorem}

  In our next section we look at some cases for which not all dif\/ferential invariants of curves are of projective type.

\subsection{The Lagrangian Grassmannian and the  manifold of reduced  pure Spinors}
These two examples are dif\/ferent, but their dif\/ferential invariants of projective type behave
similarly and so we will present them in a joint section.

\medskip

{\bfseries \itshape Lagrangian Grassmannian.} Let $G = {\rm Sp}(2n)$. If $g \in G$ then, locally
\[
g = g_{-1}g_0g_1 = \left(\begin{array}{cc}I&u\\ 0&I\end{array}\right) \left(\begin{array}{cc}\Theta&0\\ 0&\Theta^{-T}\end{array}\right) \left(\begin{array}{cc}I&0\\ S&I\end{array}\right)
\]
with $u$ and $S$ symmetric $n\times n$ matrices and $\Theta \in {\rm GL}(n)$. Again, this factorization corresponds to the splitting given by the gradation (\ref{grad}). The subgroup $H$ is
locally def\/ined by the choice $u = 0$ and a local section of the quotient can be represented by
$g_{-1}$. As usual, the action of the group is determined by the relation
\[
g \begin{pmatrix} I & u \\ 0&I\end{pmatrix} = \begin{pmatrix} I&g\cdot u\\ 0& I\end{pmatrix} h
\]
for some $h\in H$. The corresponding splitting of the algebra is given by
\[
V = V_{-1} +V_0 + V_1 = \begin{pmatrix} 0&S_1\\ 0&0\end{pmatrix}+ \begin{pmatrix} A&0\\ 0& -A^T\end{pmatrix}+\begin{pmatrix} 0&0\\ 0&S_2\end{pmatrix},
\]
where $S_1$ and $S_2$ are symmetric matrices and $A\in \gl(n)$.
The manifold $G/H$ is usually called the  \emph{Lagrangian Grassmanian} in $\R^{2n}$
and it is identif\/ied with the manifold of Lagrangian planes in $\R^{2n}$.

The following theorem describes a representation of the manifold $\K$ for curves of Lagrangian
planes in $\R^{2n}$ under the above action of ${\rm Sp}(2n)$.

\begin{theorem}[\cite{M4}] There exists a left moving frame $\rho$ along a generic curve of
Lagrangian planes such that its Serret--Frenet equations are given by
\[
K = \rho^{-1} \rho_x = \left(\begin{array}{cc} K_0 & I\\ K_1 & K_0\end{array}\right),
\]
where $K_0$ is skew-symmetric and contains all differential invariants of order~$4$, and where $K_1 = -\frac12 \s_d$. The matrix $\s_d$ is diagonal and contains in its diagonal the eigenvalues of the {\bfseries Lagrangian Schwarzian derivative} (Ovsienko {\rm \cite{Ov}})
\[
\s(u) = u_1^{-1/2}\left( u_3 - \frac32 u_2 u_1^{-1} u_2\right) (u_1^{-1/2})^T.
\]
The entries of $K_0$ and $K_1$ are functionally independent differential
invariants for curves of Lagrangian planes in $\R^{2n}$ under the action of
${\rm Sp}(2n)$. They generate
all other differential invariants. The $n$ differential invariants that
appear in $K_1$ are the invariants of projective type.
\end{theorem}

 Now we describe some of the geometric f\/lows that preserve the value $K_0=0$.
 Therefore, geometric f\/lows as below will af\/fect only invariants of projective type,
 if proper initial conditions are chosen.
 \begin{theorem}[\cite{M4}]\label{Luev}  Assume $u:J\subset \R^2 \to {\rm Sp}(2n)/H$ is a flow solution
 of
 \begin{gather}\label{Luev*}
 u_t  = \Theta^T u_1^{1/2} \, {\bf h} \, u_1^{1/2} \Theta,
 \end{gather}
 where $\Theta(x,t) \in O(n)$ is the matrix diagonalizing $\s(u)$ (i.e., $\Theta \s(u) \Theta^T = \s_d$) and where ${\bf h}$ is a symmetric matrix
 of differential invariants. Assume ${\bf h}$ is diagonal.  Then the flow
 preserves $K_0 = 0$.
 \end{theorem}

 Finally,  our next theorem gives integrable PDEs with geometric realizations as geometric f\/lows of Lagrangian planes.
 \begin{theorem}[\cite{M4}] Let $\K_1$ be the submanifold of $\K$ given by
 $K_0 = 0$. Then, the reduced bracket on $\K$ restricts to $\K_1$ to
 induce a decoupled system of $n$ second Hamiltonian structures for KdV.
 Bracket \eqref{br2} also reduces to $\K_1$ (even though
 it does not in general reduce to $\K$ for any value of $L_0$). The reduction for the choice
 \[
 L_0 = \begin{pmatrix} 0 & 0\\ I & 0\end{pmatrix} = \begin{pmatrix} 0& I\\ 0&0\end{pmatrix}^\ast \in g_1
 \]
 is a decoupled system of $n$ first KdV Hamitonian structures.

Furthermore, assume $u(t,x)$ is a flow solution of \eqref{Luev*} with ${\bf h} = \s_d$. Then \eqref{Luev*}
becomes the Lagrangian Schwarzian KdV evolution
\[
u_t = u_3 -\frac32 u_2u_1^{-1} u_2.
\]
If we choose initial conditions for which $K_0 = 0$, then the
differential invariants $\s_d$ of the flow satisfy the
 equation
 \[
 (\s_d)_t = \left({\bf D}^3 + \s_d {\bf D} + (\s_d)_x\right) {\bf h},
 \]
where ${\bf D}$ is the diagonal matrix with $\frac d{dx}$ down its diagonal.  \end{theorem}

Accordingly, if we choose ${\bf h} = \s_d$, then $\s_d$ is the solution of a \emph{decoupled system of $n$ KdV equations}
 \[
 (\s_d)_t = (\s_d)_{xxx} + 3 \s_d (\s_d)_x.
 \]

{\bfseries\itshape Reduced pure Spinors.}
A parallel description can be given for a dif\/ferent case, that of $G = O(n,n)$. In this case, if $g\in G$, locally
\[
g = g_{-1}g_0g_1 =   \begin{pmatrix} I-u&-u\\ u&I+u\end{pmatrix}\frac12 \begin{pmatrix}\Theta^{-1} + \Theta^{T}
& \Theta^{T}-\Theta^{-1}\\
\Theta^{T}-\Theta^{-1}& \Theta^{-1} + \Theta^{T}\end{pmatrix} \begin{pmatrix}I-Z&Z\\ -Z&I+Z\end{pmatrix},
\]
where $u$ and $S$ are now skew-symmetric matrices and where $\Theta\in {\rm GL}(n,\R)$.
The corresponding gradation of the algebra  is given by
$
V = V_{-1}+V_0+V_1 \in \g_{-1}\oplus\g_0\oplus\g_1$ with
\[
V_{-1}= V_{-1}(y) = \begin{pmatrix} -y&-y\\ y&y\end{pmatrix}, \qquad  V_0 = V_0(C) = \begin{pmatrix}
A&B\\ B&A\end{pmatrix}, \qquad V_1 = V_1(z) = \begin{pmatrix} -z&z\\ -z&z\end{pmatrix},
\]
$y$ and $z$ skew symmetric and $C = A+B$ given by the symmetric ($B$) and skew-symmertric
($A$) components of $C$.

Assume now that $G = O(2m,2m)$. This case has been worked out in \cite{M3}.
The homogeneous space is locally equivalent to the manifold of reduced pure spinors
in the sense of \cite{Ba}. The odd dimensional spinor case is worked out in \cite{M7}.
 Although somehow
similar it is more cumbersome to describe, so we refer the reader to \cite{M7}.

\begin{theorem}[\cite{M3}]
Let $u$ be a generic curve in $O(2m,2m)/H$. There exists a left moving frame $\rho$
such that  the left Serret--Frenet equations associated to $\rho$
are defined by
\begin{gather*}
K = V_{-1}(J)+V_0(R)+\frac18V_1(\calD),
\end{gather*}
where $J= \begin{pmatrix} 0& I_m\\ -I_m&o\end{pmatrix}$ and $R$ is of the form
\begin{gather*}
R = \begin{pmatrix} R_1& R_2\\ R_3 & -R_1^T\end{pmatrix}\in \Sp(2m)
\end{gather*}
with $R_2$ and $R_3$ symmetric, $R_1 \in \gla(m)$.  The matrix $R$ contains in the entries off the diagonals of $R_i$, $i=1,2,3$, a generating set of independent
fourth order differential invariants. The diagonals of $R_i$, $i=1,2,3$
contain a set of $3m$ independent and generating differential invariants
of order $5$ for $m > 3$ and of order $5$ and higher if $m\le 3$. The matrix $\calD$ is the skew-symmetric diagonalization of the Spinor Schwarzian
 derivative.
\end{theorem}

The  Spinor Schwarzian
 derivative  \cite{M3} is described as follows:
  if $u$ is a generic curve represented by skew-symmetric matrices, $u_1$ is non degenerate and can be brought to a normal form using $\mu\in\gla(2m)$. The matrix $\mu$ is determined up to an element
  of the symplectic group ${\rm Sp}(2m)$ (see \cite{M3}) and
 \[
\mu u_1 \mu^T = J.
\]
We def\/ine the {\it Spinor Schwarzian derivative} to be
\[
S(u) = \mu \left(u_3-\frac32 u_2 u_1^{-1}u_2\right)\mu^T,
\]
again, unique up to the action of ${\rm Sp}(2m)$. One can then prove \cite{M3} that, for
a generic curve, the Schwarzian
derivative can be diagonalize using an element of the symplectic group. That is, there exists $\theta\in {\rm Sp}(2m)$ such that
\[
\theta S(u)\theta^T = \calD = \begin{pmatrix} 0&{\bf d}\\ -{\bf d} & 0\end{pmatrix}
\]
with ${\bf d}$ diagonal. The matrix $\calD$ is the one appearing in the Serret--Frenet equations and it contains in its
entries $m$ dif\/ferential invariants {\it of projective type}.

The Spinor case seems to be dif\/ferent from others and the behavior of the Poisson brackets~(\ref{br1}) and~(\ref{br2}) is not completely understood yet. Still we do know how the invariants of projective type behave
under the analogous of the KdV Schwarzian evolution. That is described in the following theorem.
The {\it Spinor KdV Schwarzian evolution} is def\/ined by the equation
\begin{gather*}
u_t =  u_3 - \frac32 u_2 u_1^{-1} u_2.
\end{gather*}

\begin{theorem}[\cite{M3}]
Let $\rho$ is a moving frame for which normalization equations of fourth order are defined
by constants  $c_4$. Assume that, as the fourth
order invariants vanish, $[R, \widehat R] = \widehat{\widehat R}+{\it block~diagonals}$,  where
 $\widehat R$ and $\widehat{\widehat R}$ are any matrices whose only  non-zero entries are in the same position
 as the nonzero normalized entries in $R$.  Assume also that $[R,[R,\widehat R]]^d=0$ for $\widehat R$ as above, where $^d$ indicates
the diagonals in the main four blocks.

Then, if we choose initial conditions with vanishing fourth order invariants, these
remain zero under the KdV Schwarzian flow, $\D_t$ and $(R^d)_t$ decouple, and $\D$ evolves as
\[
\D_t J = \D_{xxx} J +3 \D \D_x,
\]
following a decoupled system of KdV equations.
\end{theorem}

 Although the hypothesis of this theorem seem to be very restrictive, they aren't. In fact, they are easily
 achieved when we construct moving frames in a dimension larger that~$4$, although they are more restrictive for the lower dimensions. The need for these conditions was explained in~\cite{M3}.

This is the f\/irst case where the manifold of vanishing non-projective invariants is not preserved.
Not only the vanishing of f\/ifth order invariants is not preserved, but the system {\it blows up} when
we approach the submanifold of vanishing f\/ifth order invariants. This situation (and the analogous
one for the odd dimensional case) is not well understood. It is possible that choosing normalization
equations for which $c_4$ involve derivatives of the third order dif\/ferential invariants  will produce
a better behaved moving frame. Or perhaps the f\/ifth order dif\/ferential
invariants are also of projective type and the Hamiltonian behavior is more complicated that
appears to be. The main problem understanding this case is the choice of a  moving frame that simplif\/ies
the study of the reduced evolution.
In the spinorial case such a choice is highly not trivial and we do not know whether there is no such
choice or it is just very involved. It is also not known whether or not (\ref{br1}) and (\ref{br2}) reduce to
$\K_1$. For more information see~\cite{M3} and~\cite{M7}.

The Lagrangian and Spinorial examples describe a certain type of  behavior (evolving as a~decoupled system of KdVs) of invariants of projective
type that is certainly dif\/ferent from our f\/irst example, that of $\RP^n$. Still, there is a third behavior
that is dif\/ferent from these two. These evolutions of KdV type appear in conformal manifolds. In the
conformal
case we have only two dif\/ferential invariants of projective type, and their Hamiltonian behavior is
that of a~complexly coupled system of KdVs.
That is what we describe in the next subsection.

 \subsection{Conformal case}
In this section we study the case of $G = O(p+1,q+1)$ acting on $\R^{p+q}$ as described
in \cite{Oc}. The case $q=0$ was originally studied in \cite{M5}. Using the gradation appearing in \cite{K}
we can locally factor an element of the group as $g = g_1 g_0g_{-1}$ (this is a factorization
for a right moving frame, not a left one, so it is slightly dif\/ferent from the previous examples), with $g_i \in G_i$
where
\begin{gather*}
g_{-1}(Y) = \begin{pmatrix} 1-\frac12 \nor Y & -Y_1^T & -\frac12 \nor Y &Y_2^T\\ Y_1 & I_p
& Y_1 & 0\\ \frac12 \nor Y & Y_1^T & 1+\frac12 \nor Y & -Y_2^T\\ Y_2&0&Y_2& I_q\end{pmatrix}\!,
\qquad g_0(a,b,\Theta) = \begin{pmatrix} a&0&b&0\\ 0&\Theta_{11}&0&\Theta_{12}\\
b&0&a&0\\ 0&\Theta_{21}&0&\Theta_{22}\end{pmatrix}\!,\nonumber\\
g_1(Z) = \begin{pmatrix} 1-\frac12 \nor Z & Z_1^T & \frac12 \nor Z &Z_2^T\\
-Z_1 & I_p& Z_1 & 0\\-\frac12 \nor Z & Z_1^T & 1+\frac12 \nor Z& Z_2^T\\ Z_2&0&-Z_2& I_q\end{pmatrix}.
\end{gather*}
The splitting
$Z = \begin{pmatrix}Z_1\\ Z_2\end{pmatrix}$ and $Y =
\begin{pmatrix}Y_1\\Y_2\end{pmatrix}$ is into $p$ and $q$ components, $\nor X$ is given by the f\/lat metric
of signature $(p,q)$
 and also $\Theta =
\begin{pmatrix}\Theta_{11}&\Theta_{12}\\ \Theta_{21}&\Theta_{22}
\end{pmatrix} \in O(p,q)$, $a^2-b^2 = 1$. Without loosing generality we will
assume that the f\/lat metric is given by $\nor X = X^T J X$, where $J = \begin{pmatrix} I_p&0\\ 0&-I_q\end{pmatrix}$.
The corresponding splitting in the
algebra is given by
\begin{gather*} V_{-1}(y) = \begin{pmatrix} 0 & -y_1^T & 0 & y_2^T\\ y_1 & 0
& y_1 & 0\\
0 & y_1^T & 0 & -y_2^T \\ y_2 & 0 & y_2 & 0\end{pmatrix} , \qquad V_0(\alpha, A)
= \begin{pmatrix} 0 & 0 & \alpha & 0 \\ 0 & A_{11} & 0 & A_{12} \\
\alpha & 0 & 0 & 0 \\ 0 & A_{21} & 0 & A_{22} \end{pmatrix}, \nonumber \\
V_1(z) =\begin{pmatrix} 0 & z_1^T & 0 & z_2^T \\
-z_1 & 0 & z_1 & 0 \\
0 & z_1^T & 0 & z_2^T \\ z_2 & 0 & -z_2 & 0 \end{pmatrix},
\end{gather*}
where $y = \begin{pmatrix} y_1\\ y_2\end{pmatrix}$, $z = \begin{pmatrix}z_1\\ z_2\end{pmatrix}$ are
the $p$ and $q$ components and  where $ A = (A_{ij}) \in {\mathfrak o}(p,q)$.
The algebra structure can be described as
\begin{gather*}  [ V_0(\alpha, A), V_1(z)] = V_1(JAJz + \alpha z), \qquad
[V_0(\alpha, A), V_{-1}(y)] = V_{-1}(Ay - \alpha y), \nonumber\\
[V_1(z), V_{-1}(y)] = 2 V_0\left(z^Ty,
Jzy^TJ-yz^T\right),\qquad
[V_0(\alpha, A), V_0(\beta, B)] = V_0(0, [A,B]).
\end{gather*}
With this factorization one chooses $H$ to be def\/ined by $Y=0$ and uses $G_{-1}$ as
a local section of $G/H$ (as such $Y = -u$, not $u$).

As before, the following theorem describes convenient choices of moving frames.
\begin{theorem}[\cite{M3}] There exists a left moving frame $\rho$ such that the Serret--Frenet equation for $\rho$
 is given by
$ \rho^{-1}\rho_x =  K = K_1 + K_0 + K_{-1}$ where
$K_{-1} = V_{-1}(e_1)$, $K_1 =  V_1(\inv_1e_1+\inv_2e_2)$ and $K_0 = V_0(0, \hat{K_0})$
with
\[
\hat K_0 = \begin{pmatrix} A_0&B_0\\B_0^T&0
\end{pmatrix},
\]
 and
\[
A_0 = \begin{pmatrix} 0&0&0&0&\dots&0\\ 0&0&-\inv_3&-\inv_4&
\dots&-\inv_p\\ 0&\inv_3&0&0&
\dots&0\\ \vdots& \vdots&\vdots&\vdots&\vdots&0\\
0&\inv_p&0&\dots&0&0\end{pmatrix},\qquad B_0 =
\begin{pmatrix} 0&0&\dots&0\\
-\inv_{p+1}&-\inv_{p+2}&\dots&-\inv_{p+q}\\
0&0&\dots&0\\
\vdots
&\vdots&\vdots&\vdots\\
0&0&\dots&0\end{pmatrix}.
\]
\end{theorem}

In this case there are only two generating dif\/ferential invariants of projective type, namely~$k_1$ and~$k_2$. Again, the behavior of the Poisson brackets
(\ref{br1}) and (\ref{br2})
with respect to this submanifold is spotless.
\begin{theorem}[\cite{M3}] Let $\K_1$ be the submanifold of $\K$ given by
 $K_0 = 0$. Then, the reduction of \eqref{br1} to $\K$ restricts to $\K_1$ to
 induce the second Hamiltonian structure for a complexly coupled KdV system.
 Bracket \eqref{br2} also reduced to $\K_1$  to produce the
 first Hamiltonian structure for this system.
\end{theorem}
And, again, the geometric realization for complexly coupled KdV is found.

\begin{theorem}[\cite{M3}] Assume $u:J\subset \R^2 \to O(p+1,q+1)/H$ is a solution
 of
 \[
 u_t = h_1 T + h_2 N,
 \]
 where $T$ and $N$ are conformal tangent and normal (see {\rm \cite{M3}}) and $h_1$, $h_2$ are
 any two functions of $k_1$, $k_2$ and their derivatives. Then
 the flow has a limit as $K_0 \to 0$.
 As $K_0 \to 0$, the evolution of $k_1$ and $k_2$ becomes
 \[
 \left(\begin{array}{c} k_1\\ k_2 \end{array}\right)_t=\left(\begin{array}{cc}
 -\frac12 D^3 + k_1D+Dk_1    &  k_2 D + D k_2 \\
 k_2 D + D k_2 &\frac12 D^3 - k_1D-Dk_1
 \end{array}\right) \left(\begin{array}{c} h_1\\ h_2 \end{array}\right).
 \]
 If we choose $h_1 = k_1$ and $h_2 = k_2$, then the evolution is a complexly coupled
 system of KdV equations.
  \end{theorem}

\section{Discussion}

The aim of this paper is to review some of the known evidence linking the character of  dif\/ferential invariant of curves in homogeneous spaces and the geometric realizations of integrable systems in those
manifolds. In
particular, we have described how projective geometry and geometric realizations of KdV-type
evolutions seem to be very closely related. A similar case can perhaps be made for Schr\"odinger
f\/lows, mKdV and sine-Gordon f\/lows as linked to Riemannian geometry. As we said before, several
authors \cite{A1, KQ1, KQ2, LP1, LP2, MSW, SW, TU1, TU2}, have described geometric realizations of these evolutions on manifolds that have what amounts to be   a classical {\it natural} moving frame, i.e., a frame whose derivatives of non tangential vectors have a tangential direction. This frame appears in Riemannian manifolds
and is generated by the action of the group in f\/irst order frames. Thus, one could call the invariants
they generate invariants of Riemannian-type.

 In fact, the most interesting question is how the geometry of the manifold
itself generates these geometric realizations. And, further, if a manifold
hosts a geometric realization of an integrable system of a certain type,
does that fact have any implications for the geometry of the manifold? In the case of projective
geometry, the following conjecture due to M.~Eastwood points us in this direction.

\medskip

\noindent
{\bf Conjecture.}  {\it In this type of
symmetric spaces there exists a natural projective structure along curves that
generates Hamiltonian structures of KdV type  along some flows.}

\medskip

In the conformal case $G = O(p+1, q+1)$ the two invariants of projective type are directly
connected to invariant dif\/ferential operators that appear in the work of Bailey and Eastwood (see~\cite{BE1,BE2}).
The authors def\/ined conformal circles as  solutions of a dif\/ferential equation.  The equation
def\/ines the curves together with a preferred parametrization. The parametrizations endow conformal circles with a projective structure (theirs is an explicit proof of Cartan's observation that a curve
in a conformal manifold inherits a natural projective structure, see \cite{Ca2}). We now know that the vanishing of the dif\/ferential equation in \cite{BE1} implies the vanishing of both dif\/ferential invariants of projective
type found in~\cite{M3}.  Therefore, the complexly coupled
system of KdV equations could be generated by the projective structure on conformal curves
that Cartan originally described.

Natural projective structures on curves have only been described for the cases $O(p+1,q+1)$~\cite{BE1}
and $ {\rm SL}(p+q)$~\cite{BE2}, but they do perhaps exist for $|r|$-graded parabolic manifolds.
Thus, resolving this conjecture and its generalizations would help to understand the more general situation of parabolic
manifolds. In \cite{DS} Drinfel'd and Sokolov described many evolutions of KdV-type linked to parabolic gradations
of the Lie algebra $\g$. It would be interesting to learn if parabolic manifolds (f\/lag manifolds) can be
used to generate geometric realizations for these systems.

\pdfbookmark[1]{References}{ref}
\LastPageEnding

\end{document}